\documentclass{amsart}
\usepackage{amsthm,graphics}

\numberwithin{equation}{section}
\newtheorem{theorem}[equation]{Theorem}
\newtheorem{lemma}[equation]{Lemma}
\newtheorem{proposition}[equation]{Proposition}

\newtheorem{definition}[equation]{Definition}
\newtheorem{example}[equation]{Example}
\begin{document}

\title{The Dirichlet Problem on Quadratic Surfaces}

\author{Sheldon Axler}
\address{Department of Mathematics\\
San Francisco State University\\
San Francisco, CA 94132, USA}
\email{axler@sfsu.edu}
\urladdr{www.axler.net}

\author{Pamela Gorkin}
\address{Department of Mathematics\\
Bucknell University\\
Lewisburg, PA 17837, USA}
\email{pgorkin@bucknell.edu}

\author{Karl Voss}
\address{Department of Mathematics\\
Bucknell University\\
Lewisburg, PA 17837, USA}
\email{kvoss@bucknell.edu}
\urladdr{www.facstaff.bucknell.edu/kvoss}



\thanks{The first author was supported in part by the National Science Foundation.}

\begin{abstract}
We give a fast, exact algorithm for solving Dirichlet problems with polynomial boundary functions on quadratic surfaces in $\mathbf{R}^n$ such as ellipsoids, elliptic cylinders, and paraboloids. To produce this algorithm, first we show that every polynomial in $\mathbf{R}^n$ can be uniquely written as the sum of a harmonic function and a polynomial multiple of a quadratic function, thus extending a theorem of Ernst Fischer. We then use this decomposition to reduce the Dirichlet problem to a manageable system of linear equations. The algorithm requires differentiation of the boundary function, but no integration. We also show that the polynomial solution produced by our algorithm is the unique polynomial solution, even on unbounded domains such as elliptic cylinders and paraboloids.
\end{abstract}

\maketitle

\section{Introduction}

In this paper we present a fast, exact algorithm for solving Dirichlet problems with polynomial boundary functions on a quadratic surface in $\mathbf{R}^n$  $(n \ge 2)$.
To illustrate the kind of Dirichlet problem we study, fix $b = (b_1, \ldots, b_n)\in \mathbf{R}^n$. For $x = (x_1, \ldots, x_n) \in \mathbf{R}^n$, we will write
\[
\|b x\|^2 = b_1^2 x_1^2 + \cdots + b_n^2 x_n^2.
\]
Suppose we are given a polynomial $p$ on $\mathbf{R}^n$.  We wish to find a harmonic polynomial that equals $p$ on the quadratic surface $\{x \in \mathbf{R}^n: \|b x\|^2 =  1\}$.

Even if all the $b_j$ are nonzero (so that our quadratic surface is bounded and is, in particular, an ellipsoid), computing a solution to this Dirichlet problem presents several difficulties. A standard means of expressing the solution to the Dirichlet problem for bounded domains involves the Green's function and integration. However, the Green's function of an ellipsoid does not have a known formula allowing for exact computations. An alternative approach avoids integration by employing a finite difference, finite element, or Galerkin-type scheme to approximate the solution, but again this procedure will not produce an exact solution.

If $b_1 = b_2 = \dots = b_n = 1$, then our quadratic surface is the unit sphere. In this case, a fast algorithm for finding exact solutions is presented in \cite{AXLER}. That algorithm involves differentiation but no integration. The basis of that algorithm is that any polynomial $p$ of degree $m$ on $\textbf{R}^n$ can be decomposed in the form
\[
             p=h+(\|x\|^2-1)f,
\]
where $h$ is a harmonic polynomial of degree at most $m$ and $f$ is a polynomial of degree at most $m-2$. Because $h$ is harmonic and equals $p$ on the unit sphere, it is the solution to our Dirichlet problem. The algorithm presented in \cite{AXLER} shows how the polynomials $h$ and $f$ in the decomposition above can be computed via differentiation from $p$. Unfortunately these techniques work only on spheres and so do not provide an algorithm for nonspherical ellipsoids or other quadratic surfaces.

In this paper we solve the Dirichlet problem discussed above, getting solutions for ellipsoids as well as for elliptic cylinders (for example, $\{x \in \mathbf{R}^3: x_1^2+2 x_2^2=1\}$) and paraboloids (for example, $\{x \in \mathbf{R}^3: x_3=x_1+x_2^2\}$). We will begin by extending the decomposition above
to a collection of quadratic surfaces. We then use this decomposition to produce a system of linear equations whose solution will give an exact solution to our Dirichlet problem. We will show how this system of linear equations has a structure allowing it to be reduced to smaller systems of linear equations, thus producing a fast algorithm. The algorithm requires differentiation of the boundary function, but no integration.

Before we turn to these matters, we need to present some background, much of which appears in \cite{SHAPIRO}. A \textit{multi-index} is an $n$-tuple $\alpha = (\alpha_1, \dots, \alpha_n)$ of nonnegative integers. 
The \textit{order} of $\alpha$, denoted $|\alpha|$, is defined by
\[
|\alpha| = \alpha_1 + \dots + \alpha_n.
\]
We let $x^\alpha$ denote  the monomial  ${x_1}^{\alpha_1} \dots {x_n}^{\alpha_n}$ and $D^\alpha$ denote the differential operator
${D_1}^{\alpha_1} \dots {D_n}^{\alpha_n}$, where $D_j$ denotes differentiation with respect to $x_j$. If $q$ is a polynomial on $\mathbf{R}^n$ given by $q(x) = \sum_\alpha c_\alpha x^\alpha$, then $q(D)$ is the differential operator defined by $q(D) = \sum_\alpha c_\alpha D^\alpha$.
A polynomial is called \textit{homogeneous} of degree $m$ if it is a linear combination of monomials of degree~$m$.

Ernst Fischer proved that given a homogeneous polynomial $q$ on $\mathbf{R}^n$, every homogeneous polynomial $p$ of degree $m$ can be decomposed uniquely as $p=h+ qf$, where $h$ is a homogeneous polynomial of degree $m$ satisfying $q(D)h=0$ and $f$ is a homogeneous polynomial of suitable degree.

In \cite{SHAPIRO} the subject of more general decompositions is discussed.  Given two polynomials $g$ and $q$ on $\mathbf{R}^n$, the relevant question is whether an arbitrary polynomial $p$ can be decomposed as $p=h + qf$, where $h$ and $f$ are polynomials, with $h$ satisfying $g(D)h=0$. Shapiro refers to a pair $(g, q)$ with this property as a \textit{generalized \mbox{Fischer} pair}. He asks: Which $(g,q)$ form generalized Fischer pairs? Note that if $g(x) = \|x\|^2$, then $g(D)$ is the Laplacian and so this decomposition would require $h$ to be harmonic. We will provide examples of a robust class of quadratic polynomials $q$ that form a generalized Fischer pair with $g(x)=\|x\|^2$, and we give explicit
examples of the decomposition via our algorithm.

Some of the surfaces that we consider are unbounded (for example, the elliptic cylinders  and paraboloids mentioned above). Thus unique solutions to Dirichlet problems on these surfaces, even in the class of polynomials, are neither automatic nor expected.  For example, the set of harmonic polynomials that vanish on the hyperplane $\{x \in \mathbf{R}^n: x_n=0\}$ is not trivial. However, we will show that for the quadratic surfaces we consider, polynomial solutions to the Dirichlet problem with polynomial boundary functions are unique.

The paper is organized as follows:  In Section~\ref{S:Fischer'slemma} we begin by presenting Fischer's lemma and the corresponding decomposition theorem. We then extend these results to cover a wider class of generalized Fischer pairs. These generalized Fischer pairs are then used to solve the Dirichlet problem. We prove that the polynomial produced by this technique is the unique polynomial solving the Dirichlet problem, even when our quadratic surface is unbounded. In Section~\ref{S:algorithm} we describe a fast algorithm for computing the solution to the Dirichlet problem promised by the results in Section~\ref{S:Fischer'slemma}. In Section~\ref{S:examples} we present some examples, computed using our algorithm, of solutions to Dirichlet problems on ellipsoids. The Appendix contains a differentiation formula needed by our algorithm.

\section{Fischer's Lemma and the Dirichlet Problem} \label{S:Fischer'slemma}

In this section, we state and prove our generalization of Fischer's results.
Then we prove the decomposition theorem that solves the Dirichlet problem. Then we show that even when our quadratic surfaces are unbounded, the solution given by our decomposition theorem is the only polynomial solution.

We begin by stating Fischer's results, which are nicely restated and proved in \cite{BAKER}; also see \cite{KhSh}. Fix an integer $n \ge 2$. We will always use $m$ to denote a nonnegative integer. Let $\mathcal{P}_m$ denote the vector space of polynomials (with real coefficients) of degree at most $m$ on~$\mathbf{R}^n$. For convenience, we define $\mathcal{P}_k$ to be $\{0\}$ for $k < 0$. As usual, $\Delta$ denotes the Laplacian.

\begin{lemma}[Fischer's Lemma]
Suppose $b = (b_1, \dots, b_n)$, where each $b_j \ne 0$.
Define $L \colon \mathcal{P}_m \to \mathcal{P}_m$ by
\[
L(f) = \Delta\bigl((\|bx\|^2 - 1) f\bigr).
\]
Then $L$ is a linear, degree-preserving, bijection of $\mathcal{P}_m$ onto itself.
\label{fischer1}
\end{lemma}

Fischer's Lemma leads to Fischer's Decomposition Theorem, which gives a solution to the Dirichlet problem for ellipsoids.

\begin{theorem}[Fischer's Decomposition Theorem]   \label{baker'stheorem}
Suppose $b = (b_1, \dots, b_n)$, where each $b_j \ne 0$.
Let $p\in \mathcal{P}_m$. Then there exists a unique harmonic polynomial
$h \in \mathcal{P}_m$ such that
\[
          \label{decomposition}
                 p = h + (\|bx\|^2-1) f
\]
for some $f \in \mathcal{P}_{m - 2}$.

Let $E = \{x \in \mathbf{R}^n : \|b x\|^2 < 1\}$. Then $h$ is the unique continuous function on $\bar{E}$ that is harmonic on $E$ and equals $p$ on the ellipsoid $\partial E$.
\end{theorem}

We will need the following generalization of Corollary 5.3 of \cite{AXLERBOOK}.

\begin{lemma}     \label{L:noellharm}
If $b = (b_1, \dots, b_n)$, where each $b_j \ne 0$, then no nonzero polynomial multiple of $\|bx\|^2$ is harmonic.
\end{lemma}

\begin{proof}
Suppose, to the contrary, that $p$ is a nonzero polynomial of degree $m$ such that $\|bx\|^2p$ is harmonic. Let $E$ be as in the theorem above, which states that 
there exists a harmonic polynomial $h \in \mathcal{P}_{m}$ such that $h$
equals $p$ on the ellipsoid $\partial E$.
But $\|bx\|^2p$ is also a harmonic polynomial that equals $p$ on $\partial E$. 
Because $\|bx\|^2 p$ has degree $m+2$, the polynomials $\|bx\|^2 p$ and $h$ cannot be equal. But this contradicts the uniqueness of solutions to the Dirichlet problem on bounded domains.
\end{proof}

In order to generate new generalized Fischer pairs, let $q$ be a quadratic (degree~$2$) polynomial on $\mathbf{R}^n$. We want to look at the map on $\mathcal{P}_m$ defined by
\[
            L(f) = \Delta (qf).
\]
Our goal is to identify choices of $q$ for which $L$ is a bijection of $\mathcal{P}_m$ onto $\mathcal{P}_m$. This leads us to the definition of a nonhyperbolic quadratic.

\begin{definition}   \label{D:nonhyperbolic}
A \textit{nonhyperbolic} quadratic is a polynomial $q$ on $\mathbf{R}^n$ of the form
\[
\sum_{j = 1}^n b_j^2 x_j^2 + \sum_{j=1}^n c_jx_j + d,
\]
where at least one $b_j \ne 0$.
\end{definition}

Note that the next theorem, which gives the desired bijectivity, implies that no nonzero polynomial multiple of a nonhyperbolic quadratic is harmonic, generalizing Lemma~\ref{L:noellharm}. This result does not hold for arbitrary quadratic polynomials. In fact, even for a nonharmonic quadratic polynomial, a nonzero polynomial multiple might be harmonic. For example, $x_1^2 - 3x_2^2$ is not harmonic, but $(x_1^2 - 3x_2^2)x_1$ is harmonic.

\begin{theorem}  \label{T:trivialkernel}
Let $q$ be a nonhyperbolic quadratic.  Define $L\colon \mathcal{P}_m \to \mathcal{P}_m$ by
\[
L(f) = \Delta(q f).
\]
Then $L$ is a linear bijection of $\mathcal{P}_m$ onto $\mathcal{P}_m$.
\end{theorem}

\begin{proof}
Clearly $L$ is a linear map of $\mathcal{P}_m$ into $\mathcal{P}_m$. Since $\mathcal{P}_m$ is finite dimensional,
we need only show that $L$ is injective. So suppose that this is not true. Then there exists $f \in \mathcal{P}_m$, $f \ne 0$, such that $L(f) = 0$. Without loss of generality, we may suppose that $f$ is of degree~$m$ (otherwise, replace $m$ by a lower integer).
We can write $f = f_m + f_{m - 1}$, where $f_m$ is homogeneous of degree~$m$ and $f_{m - 1} \in \mathcal{P}_{m - 1}$. Write $q = q_2 + q_1$, where $q_2$ is homogeneous degree $2$ (so $q_2(x) =
\sum_{j = 1}^n b_j^2 x_j^2$ in the notation above) and $q_1\in \mathcal{P}_1$. Because $\Delta(q f) = 0$, we know that
\[
\Delta(q_2 f_m + q_2 f_{m-1} + q_1 f) = 0.
\]
Thus
\[
\Delta(q_2 f_m) = 0,
\]
because the other terms have lower degrees. We will now show that this implies
that $f_m = 0$, which is a contradiction.

Note that we have reduced our theorem to the case where $q = q_2$. Reordering the variables (if some of the $b_j = 0$), we see that it suffices to prove our theorem in the case when $q(x) = \sum_{j = 1}^r b_j^2 x_j^2$, where $1 \le r \le n$ and $b_1, \dots, b_r$ are all nonzero. To simplify notation, we will also replace $f_m$ in the previous paragraph with $f$. So again we have the assumption that $\Delta(qf) = 0$ and we want to prove that $f = 0$, but now we have a special form for $q$.

If $r = n$, our desired conclusion that $f = 0$ follows from Lemma~\ref{L:noellharm}. So suppose $r < n$. Let $k$ denote the degree of $f$ thought of as a polynomial in $x_{r + 1}, \ldots, x_n$ (temporarily think of $x_1, \dots, x_r$ as constants to define $k$). Write
\begin{equation} \label{decomp26}
         f = p + g,
\end{equation}
where $p$ is the part of $f$ that is homogeneous of degree $k$ in the variables
$x_{r + 1}, \ldots, x_n$ and $g$ is the remaining part of $f$, consisting of lower degree terms in the variables $x_{r + 1}, \ldots, x_n$.

Using the product formula for the Laplacian, which states that
\begin{equation} \label{prod}
\Delta(q p) =  p \Delta q +q \Delta p 
+ 2 \nabla q \cdot \nabla p,
\end{equation}
we obtain from (\ref{decomp26}) the equation
\[
\Delta(q f) = p \Delta q + q \Delta p
+ 2 \nabla q \cdot \nabla p + \Delta(qg).
\]
Let's consider the degree, as a polynomial in $x_{r + 1}, \ldots, x_n$, of each term on the right side of this equation. Because $\Delta q$ is a positive constant, this degree is $k$ for the first term. Because $q$ is independent of the variables $x_{r + 1}, \ldots, x_n$, the second and third terms have degree (as a function of $x_{r + 1}, \ldots, x_n$) less than $k$. Because the degree of $g$ (as a function of
$x_{r + 1}, \ldots, x_n$) is less than $k$ and $q$ is independent of
$x_{r + 1}, \ldots, x_n$,  the degree of the fourth term is less than $k$.

Thus the only part of the right side of the equation above with degree $k$ (as a function of $x_{r + 1}, \ldots, x_n$) is the first term, $p \Delta q$. The left side of the equation is $0$, so $p \Delta q = 0$. Hence $p = 0$. But $p$ was the part of $f$ of highest degree in $x_{r + 1}, \ldots, x_n$. Hence $f$ is independent of $x_{r + 1}, \ldots, x_n$. Thus we can think of the equation $\Delta (qf) = 0$ as taking place in $\mathbf{R}^r$. Lemma~\ref{L:noellharm} now implies that $f = 0$, as desired.
\end{proof}

We now apply Theorem~\ref{T:trivialkernel} to obtain the general decomposition  theorem.

\begin{theorem}  \label{T:decomposition}
Suppose $p \in \mathcal{P}_m$.  Let $q$ be a nonhyperbolic
quadratic. Then there exists a unique harmonic polynomial $h \in \mathcal{P}_m$ such that
\[
            p= h+ qf
\]
for some $f \in \mathcal{P}_{m-2}$.
\end{theorem}

\begin{proof}
Note that $\Delta p \in \mathcal{P}_{m-2}$. Thus by Theorem~\ref{T:trivialkernel}, there exists $f \in \mathcal{P}_{m-2}$
such that $\Delta (qf) = \Delta p$.  Let $h=p-qf$. Then $h$ is harmonic polynomial in $\mathcal{P}_m$ and $p = h + qf$, as desired.

To prove the uniqueness part of this theorem, suppose also that $\tilde{h}$ is a harmonic polynomial in $\mathcal{P}_m$ and that $p = \tilde{h} + q\tilde{f}$ for some $\tilde{f} \in \mathcal{P}_{m-2}$. Then
\[
h - \tilde{h} = q(\tilde{f} - f).
\]
The left side of the equation above is harmonic, and hence $\Delta \bigl(q(\tilde{f} - f)\bigr) = 0$. Theorem~\ref{T:trivialkernel} now implies that $\tilde{f} - f = 0$, which implies that $\tilde{h} = h$, as desired.
\end{proof}

If $p \in \mathcal{P}_m$ and $q$ is a  nonhyperbolic
quadratic, we can consider the following Dirichlet problem: find a harmonic polynomial $h \in \mathcal{P}_m$ such that $h$ equals $p$ on the set $\{x \in \mathbf{R}^n : q(x) = 0 \}$. Clearly the $h$ produced by the theorem above solves this Dirichlet problem. Of course, $\{x \in \mathbf{R}^n : q(x) = 0 \}$ could be the empty set or a single point. Because $q(x) \to \infty$ as $|x| \to \infty$, the existence of a point in $\mathbf{R}^n$ where $q$ is negative is a convenient condition to ensure that $\{x \in \mathbf{R}^n : q(x) = 0 \}$ is a nondegenerate quadratic surface.
For example, using the notation of Definition~\ref{D:nonhyperbolic}, $\{x \in \mathbf{R}^n : q(x) = 0 \}$ will be a nondegenerate quadratic surface if
\[
d < \sum_{\{j:b_j \ne 0\}} \frac{c_j^2}{4b_j^2}
\]
or if $c_j \ne 0$ for some $j$ with $b_j = 0$.

We now turn to the question of whether the polynomial $h$ produced by Theorem~\ref{T:decomposition} is the unique polynomial solution to the Dirichlet problem discussed in the paragraph above. The following lemma will help us answer this uniqueness question. If we were working in $\mathbf{C}^n$ instead of $\mathbf{R}^n$, then Hilbert's Nullstellensatz could be used to provide information about when a polynomial $h$ vanishing on the zero set of another polynomial $q$ is a polynomial multiple of $q$. A theorem called the Real Nullstellensatz (see, for example, \cite{Scharlau}, Chapter~3, Theorem~3.3) provides some information about polynomials on $\mathbf{R}^n$ vanishing on the zero set of another polynomial. However, we do not see how the Real Nullstellensatz can be used to prove the lemma below, so we have provided a proof without using such machinery.

\begin{lemma} \label{unique}
Suppose $q$ is a nonhyperbolic quadratic that is negative at some point of $\mathbf{R}^n$. If $h$ is a polynomial on $\mathbf{R}^n$ such that $h(x) = 0$ whenever $q(x) = 0$, then $h$ is a polynomial multiple of $q$.
\end{lemma}

\begin{proof}
We will prove this lemma by induction on the dimension $n$.

To get started, suppose $n = 1$ and that $q(x) = b^2 x^2 + cx + d$, where $b \ne 0$, is a nonhyperbolic quadratic that is negative at some point of $\mathbf{R}$. Suppose $h$ is a polynomial on $\mathbf{R}$ such that $h(x) = 0$ whenever $q(x) = 0$. Because $q$ is negative at some point of $\mathbf{R}$ and $q(x) \to \infty$ as $x \to \infty$, we see that $q$ has precisely two distinct zeros. The polynomial $h$ vanishes on both these zeros, and thus $h$ is a polynomial multiple of the quadratic polynomial $q$, as desired.

Now suppose that the lemma holds in dimension $n-1$. Let $q$ be a nonhyperbolic quadratic that is negative at some point of $\mathbf{R}^n$. Relabelling coordinates, if necessary, we can assume that
\[
q(x) =
\sum_{j = 1}^n b_j^2 x_j^2 + \sum_{j=1}^n c_jx_j + d,
\]
where $b_j \ne 0$ for some $j \in \{1, \dots, n-1\}$. Let $y$ denote a typical point of $\mathbf{R}^{n-1}$ and let $z$ denote a typical point of $\mathbf{R}$; thus $(y,z)$ denotes a typical point of $\mathbf{R}^n$. 

Suppose $h$ is a polynomial on $\mathbf{R}^n$ of degree $m$ such that $h(x) = 0$ whenever $q(x) = 0$. We need to show that $h$ is a polynomial multiple of $q$.

For $z \in \mathbf{R}$ such that $q(y, z)$ is negative for some $y \in \mathbf{R}^{n-1}$, define a nonhyperbolic quadratic $q_z$ on $\mathbf{R}^{n-1}$ by
\[
q_z(y) = q(y, z),
\]
and define a polynomial $h_z$ on $\mathbf{R}^{n-1}$ of degree at most $m$ by
\[
h_z(y) = h(y, z).
\]
Then $h_z(y) = 0$ whenever $q_z(y) = 0$, and thus by our induction hypothesis there is a polynomial $f_z$ on $\mathbf{R}^{n-1}$ such that
\begin{equation} \label{uniq1}
h_z(y) = f_z(y) q_z(y)
\end{equation}
for all $y \in \mathbf{R}^{n-1}$. Clearly $f_z$ has degree at most $m-2$.

For $y \in \mathbf{R}^{n-1}$ such that $q(y, z)$ is negative for some $z \in \mathbf{R}$, define a  polynomial $q^y$ on $\mathbf{R}$ by
\[
q^y(z) = q(y, z),
\]
and define a polynomial $h^y$ on $\mathbf{R}$ of degree at most $m$ by
\[
h^y(z) = h(y, z).
\]
Then $h^y(z) = 0$ whenever $q^y(z) = 0$. We claim that there is a polynomial $g^y$ on $\mathbf{R}$ such that
\begin{equation} \label{uniq2}
h^y(z) = g^y(z) q^y(z)
\end{equation}
for all $z \in \mathbf{R}$. If $b_n \ne 0$, then $q^y$ is a nonhyperbolic quadratic on $\mathbf{R}$ and the claim follows from the dimension~$1$ case that was proved at the beginning of this proof. If $b_n = 0$ but $c_n \ne 0$, then $q^y$ is a polynomial on $\mathbf{R}$ of degree 1 and the claim follows easily. Finally, if $b_n = 0$ and $c_n = 0$, then $q^y$ is a negative constant, in which case the claim is trivially true. In any case, we see that $g^y$ has degree at most $m$.

Let $\Omega = \{(y, z) \in \mathbf{R}^n: q(y, z) < 0\}$. Combining (\ref{uniq1}) and (\ref{uniq2}), we see that
\[
h(y,z) = f_z(y) q(y, z) = g^y(z) q(y, z)
\]
for all $(y, z) \in \Omega$. Thus we can define a function $p$ on $\Omega$ by
\[
p(y, z) = f_z(y) = g^y(z).
\]
Because
\[
p(y, z) = \frac{h(y,z)}{q(y,z)}
\]
for all $(y, z) \in \Omega$, we see that $p$ is real-analytic on $\Omega$.

Suppose $\alpha = (\alpha_1, \dots, \alpha_n)$ is a multi-index of order greater than $2m$. Then either
\[
\alpha_1 + \dots + \alpha_{n-1} > m-1 \quad \text{or} \quad \alpha_n > m+1. 
\]
This implies that $D^\alpha p(y,z) = 0$ for all $(y,z) \in \Omega$. Because all sufficiently high-order partial derivatives of $p$ equal $0$ on the open set $\Omega$, we conclude that $p$ is a polynomial on $\Omega$. Hence we can think of $p$ as a polynomial defined everywhere on~$\mathbf{R}^n$.

Finally, because $h(y, z) = p(y, z) q(y, z)$ for all $(y, z) \in \Omega$, and because polynomials that agree on a nonempty open subset of $\mathbf{R}^n$ must agree everywhere, we have $h = p q$. Thus $h$ is a polynomial multiple of $q$, as desired.
\end{proof}

Now we can combine the previous lemma and Theorem~\ref{T:trivialkernel} to prove the desired uniqueness result. Of course  in the ellipsoidal case (where each $b_j \ne 0$, in the notation of Definition~\ref{D:nonhyperbolic}) uniqueness follows easily from the boundedness of the surface, but we want to consider also elliptic cylinders and paraboloids. Note that the uniqueness result in the theorem below fails on some nondegenerate quadratic surfaces, so the hypothesis that $q$ is nonhyperbolic cannot be deleted. For example, on the quadratic surface defined by $\{x \in \mathbf{R}^n: x_1^2 - 3x_2^2 - 1 = 0\}$, any solution to a Dirichlet problem can be used to produce another solution by adding to it the harmonic polynomial $(x_1^2 - 3x_2^2 - 1)x_1$, which vanishes on the quadratic surface in question.

To obtain uniqueness results on half-spaces, even in the class of polynomial solutions, a growth condition on the solutions is needed (see~\cite{SIEGELTALVILA}). However, the theorem below shows that we have unique polynomial solutions on our quadratic surfaces without the requirement of a growth condition.

\begin{theorem}
Suppose $q$ is a nonhyperbolic quadratic that is negative at some point of $\mathbf{R}^n$. If $p \in \mathcal{P}_m$, then there is a unique harmonic polynomial $h$ that equals $p$ on $\{x \in \mathbf{R}^n: q(x) = 0\}$. Furthermore,
\[
h = p - qf
\]
for some $f \in \mathcal{P}_{m-2}$.
\end{theorem}

\begin{proof}
Take $h$ and $f$ as in Theorem~\ref{T:decomposition}. Then $h$ is a harmonic polynomial that equals $p$ on $\{x \in \mathbf{R}^n: q(x) = 0\}$; furthermore
$f \in \mathcal{P}_{m-2}$ and $h = p - qf$.

To prove uniqueness, suppose $\tilde{h}$ is also a harmonic polynomial that equals $p$ on $\{x \in \mathbf{R}^n: q(x) = 0\}$. Then $\tilde{h} - h$ equals $0$ on $\{x \in \mathbf{R}^n: q(x) = 0\}$. By Lemma~\ref{unique}, $\tilde{h}- h$ is a polynomial multiple of $q$. However, Theorem~\ref{T:trivialkernel} implies that no nonzero polynomial multiple of $q$ is harmonic. Thus $\tilde{h}- h = 0$ and hence $\tilde{h} = h$, completing the proof of uniqueness.
\end{proof}

\section{Algorithm} \label{S:algorithm}
 
In this section, we will turn the results of the previous section into a
computationally useful algorithm for solving the Dirichlet problem. Every polynomial is a sum of homogeneous polynomials, and thus to solve Dirichlet problems with polynomial boundary functions it suffices to do so for homogeneous polynomials.

Let $\mathcal{H}_m$ denote the vector space of polynomials on $\mathbf{R}^n$ that are homogeneuos of degree $m$  (the polynomial $0$ is homogeneous of every degree, so $0 \in  \mathcal{H}_m$).
Suppose $q$ is a nonhyperbolic quadratic and $p\in \mathcal{H}_{m+2}$.  We use the
decomposition in Theorem~\ref{T:decomposition} to write $p = h + q f$, where $h$ is a harmonic polynomial in $\mathcal{P}_{m+2}$ and $f \in \mathcal{P}_m$.
Breaking each polynomial into its homogeneous components, we obtain
\[
         p= \sum_{j = 0}^{m+2} h_j+(q_2+q_1+q_0)\sum_{j = 0}^m f_j,
\]
where each $h_j$ is harmonic.  We can break the equation above
into homogeneuous equations by degree to obtain the system
   \begin{align}
         p& = h_{m+2}+q_2 f_m, \notag\\
         -q_1 f_m& = h_{m+1}+q_2 f_{m-1}, \notag\\
         -q_0 f_m-q_1 f_{m-1}& = h_m+q_2f_{m-2}, \notag\\
         &\vdots \label{homogen}\\
         -q_0 f_2-q_1 f_1& = h_2+q_2f_0, \notag\\
         -q_0 f_1-q_1 f_0 & =  h_1, \notag\\
         -q_0f_0& =  h_0 .\notag
   \end{align}
Here $p$ and $q_2, q_1, q_0$ are known and we need to compute $h_{m+2}, \dots, h_0$ and $f_m, \dots, f_0$. Our plan of attack is to use the first equation to find $f_m$, which will then give us both $h_{m+2}$ and the left side of the next equation.  We can repeat the procedure with the second equation to find $f_{m-1}$, which will give us both $h_{m+1}$ and the left side of the next equation.  We continue this process until we have found $h_j$ for each $j$. Then $h=\sum _{j=0}^{m+2} h_j$ is the harmonic function that agrees with $p$ on $\{x \in \mathbf{R}^n: q(x) = 0\}$.

Thus, we turn our attention to equations of the form
    \begin{equation} \label{E:general}
          p=h+q_2f
    \end{equation}
where $p$ is a known homogeneous polynomial, $q_2$ is the highest degree part of a known nonhyperbolic quadratic, $h$ is an unknown harmonic polynomial homogeneous of degree $\deg(p)$, and $f$ is an unknown polynomial homogeneous of degree $\deg(p) - 2$. Note that all the equations in the system above are of this form (except the last two equations, which are trivial to solve for $h_1$ and~$h_0$). Hence an algorithm for finding the solution to (\ref{E:general}) will give us an algorithm for solving our Dirichlet problem.

To solve equation (\ref{E:general}), eliminate $h$ by taking the Laplacian of both sides, getting
\begin{equation}\label{3}
\Delta p = \Delta(q_2f).
\end{equation}
Now our problem has been reduced to finding $f \in \mathcal{H}_m$ satisfying the equation above, where $p \in \mathcal{H}_{m+2}$ and $q_2(x) = b_1^2 x_1^2 + \dots + b_n^2 x_n^2$ are known, with at least one $b_j \ne 0$. Note that Theorem~\ref{T:trivialkernel} implies that (\ref{3}) has a unique solution $f \in \mathcal{H}_m$.

The algorithm we introduce for finding $f$ involves repeated differentiation of~(\ref{3}). Because $f \in \mathcal{H}_m$, we will know $f$ once we know the constants
$D^\alpha f$ for every multi-index $\alpha$ of order~$m$.
We need one more piece of notation before finding the constants $D^\alpha f$. For $1 \le j \le n$, let $e_j$ denote the multi-index whose $j^{\text{th}}$-coordinate equals $1$ and whose other coordinates equal $0$.

For $\alpha$ a multi-index of order $m$, apply the differential operator $D^\alpha$ to both sides of (\ref{3}), getting
\begin{align}
D^\alpha(\Delta p) &= D^\alpha\bigl(\Delta(q_2f)\bigr) \label{dalpha} \\
&= \Delta \bigl(D^\alpha (q_2 f)\bigr) \notag \\
&= \Delta \Bigl(q_2 D^\alpha f + \sum
               _{j=1}^n 2 \alpha_j b_j^2 x_j D^{\alpha - e_j }
               f+
\sum _{j=1}^n  \alpha_j (\alpha_j - 1) b_j^2 D^{\alpha - 2e_j} f \Bigr), \notag
\end{align}
where the last equality comes from Proposition~\ref{theorem:1} in our Appendix and the explicit form of $q_2$. Because $D^\alpha f$ is a constant, the Laplacian of the first term in parentheses above equals $2 \|b\|^2 D^\alpha f$. Because $D^{\alpha-e_j}f$ has degree 1, the Laplacian of the second term in parentheses above can be easily computed using the product formula for the Laplacian (\ref{prod}), and the last equation becomes
\begin{equation} \label{sys}
           D^\alpha(\Delta p)
            = \bigl(2 \|b\|^2 + 4 \sum_{j = 1}^n \alpha_j b_j^2\bigr) D^\alpha f+
            \sum_{j=1}^n  \bigl[ \alpha_j (\alpha_j -1) b_j^2
               \sum_{k=1}^n (D^{\alpha +2e_k - 2e_j} f)\bigr].
\end{equation}

Note that multi-index $\alpha + 2e_k - 2e_j$ in the equation above has order $m$. 
Thus as $\alpha$ ranges over all multi-indices of order $m$, (\ref{sys}) gives us a system of linear equations in the unknowns $\{D^\alpha f\}$. This system of equations can be solved using Gaussian elimination, giving us $f$ and thus solving our Dirichlet problem.

Any solution to the system of equations (\ref{sys}) in the unknowns $\{D^\alpha f\}$ gives a function $f \in \mathcal{H}_m$ satisfying (\ref{dalpha}) for all multi-indices $\alpha$ of order $m$, which implies that $f$ satisfies $\Delta p= \Delta(q_2f)$. However, we already know from Theorem~\ref{T:trivialkernel} that there is a unique $f \in \mathcal{H}_m$ satisfying $\Delta p= \Delta(q_2f)$. Thus the system of equations (\ref{sys}) in the unknowns $\{D^\alpha f\}$ has a unique solution.

The number of operations needed to compute the solution to a system of $s$ linear equations in $s$ unknowns is on the order of $(2/3)s^3$. Thus we can expect that solving the system (\ref{sys}) will take on the order of
\begin{equation} \label{big}
\frac{2 (m+n-1)!^3}{3 m!^3 (n-1)!^3}
\end{equation}
operations (the formula for the number of multi-indices of order $m$ can be found, for example, on page~78 of \cite{AXLERBOOK}). However, a careful look at the system (\ref{sys}) shows that we can do much better. 

For a fixed multi-index $\alpha$, the multi-indices that appear in (\ref{sys}) all have the form $\alpha + 2e_k - 2e_j$. This leads us to define an equivalence relation on the set of multi-indices of order $m$ by declaring that two multi-indices $\alpha$ and $\beta$ are equivalent if $\alpha_i \equiv \beta_i \pmod{2}$ for each $i = 1, \dots, n$. This equivalence relation breaks the set of multi-indices of order $m$ into equivalence classes, and the system of equations (\ref{sys}) breaks into corresponding systems of equations. Hence instead of solving one large system of equations, we can solve several smaller systems of equations, which leads to considerable computational savings, as we will soon see.

If $m$ is larger than $n$, then we have $2^{n-1}$ equivalence classes, corresponding to a choice of even or odd entry in the first $n-1$ coordinates of a multi-index (the parity of the last coordinate is forced by the condition that the coordinates add up to $m$). As $m$ increases, the number of elements in each of these $2^{n-1}$ equivalence classes divided by the total number of multi-indices of order $m$ approaches $1/2^{n-1}$. Thus instead of solving one large system of $\binom{m+n-1}{m}$ equations in the same number of variables, we can solve $2^{n-1}$ systems of equations, each system containing on the order of $\binom{m+n-1}{m}/2^{n-1}$ equations and variables. This computation requires on the order of
\begin{equation} \label{small}
2^{n-1} \frac{2 (m+n-1)!^3}{3 (2^{n-1})^3 m!^3 (n-1)!^3}
\end{equation}
operations.

The ratio of (\ref{big}) to (\ref{small}) is $2^{2n-2}$. Thus breaking our system of equations into smaller systems of equations using our equivalence relation reduces the number of operations by a factor of $2^{2n-2}$. For example, if $n = 6$, then this technique should reduce computation time by a factor of over 1000. This savings is needed even for moderate size $m$ and $n$, because the system of equations $(\ref{sys})$ grows large rapidly.

Usually we can do even better than reducing computations by a factor of $2^{2n-2}$. Suppose, for example, that $p(x) = x_1^{20} x_2^7$. Then $\Delta p = 42 x_1^{20}x_2^5 + 380 x_1^{18}x_2^7$. Note that the left side of (\ref{sys}) equals 0 for all multi-indices of order 25 except $(20, 5, 0, \dots, 0)$ and $(18, 7, 0, \dots, 0)$. 
In other words, in only one of our equivalence classes (the equivalence class consisting of those multi-indices of order 25 whose first coordinate is even, second coordinate is odd, and coordinates 3 through $n$ are even) is the left side of (\ref{sys}) anything other than 0. When the left side of every equation in an equivalence class is 0, there is no need to perform Gaussian elimination to solve the system of equations in that equivalence class, because obviously all the unknowns equal 0 (recall that the system (\ref{sys}) has a unique solution). Thus in the example at hand, instead of solving $2^{n-1}$ smaller systems of equations, we need only solve one smaller system of equations.

As can be seen from the reasoning in the previous paragraph, if $p$ is a monomial then our computation time is reduced, through the use of our equivalence classes, by another factor of $2^{n-1}$. Thus if $p$ is a monomial, then computation time through the use of equivalence classes is reduced by a factor of $2^{3n-3}$. For example, if $n = 6$, then this technique should reduce computation time by a factor of over 32,000.

The full strength of the reduction discussed above by a factor of $2^{3n-3}$, as opposed to a reduction by a factor of $2^{2n-2}$, holds only for ellipsoids and elliptic cylinders. For ellipsoids and elliptic cylinders we can assume (perhaps after a translation) that each $c_j$ in Definition \ref{D:nonhyperbolic} equals $0$, which gives $q_1 = 0$ in the system (\ref{homogen}). However, if $q_1 \ne 0$ then multiplication by $q_1$, when solving the second and successive equations in the system (\ref{homogen}), can lead to nonzero left sides in equations in the system (\ref{sys}) other than the equations in the equivalence class corresponding to the exponents in the monomial~$p$. 

\section{Examples}  \label{S:examples}

The algorithm outlined in the previous section has been implemented by the authors in \textit{Mathematica} and in \textit{MATLAB}. The \textit{Mathematica} version produces exact solutions, in considerably less time than we believed possible when we started this project. Even with simple boundary polynomials in low dimensions, solutions to Dirichlet problems on nonhyperbolic quadratic surfaces tend to involve fractions with large numerators and denominators, as we will see in the examples presented below. Our \textit{Mathematica} software can use floating point arithmetic to produce only decimal approximations to the exact solutions, with even faster times than when working with exact rational arithmetic. Our \textit{MATLAB} implementation of the algorithm works only in floating point arithmetic, again quickly producing decimal approximations to the exact solutions.

Our \textit{Mathematica} implementation of the algorithm is available on the first author's web site; our \textit{MATLAB} implementation of the algorithm is available on the third author's web site. The appropriate web addresses are listed at the end of this paper; within those web sites look for information about this paper to find the software. Although our software is available without charge, its use requires \textit{Mathematica} or \textit{MATLAB}. The examples presented below were generated by our \textit{Mathematica} implementation of the algorithm. 

Our first example will be in dimension $3$ with a boundary function of degree~$7$. Here, as elsewhere throughout this section, all fractions are given in reduced form.

\begin{example} \label{ex1}
Suppose $p(x_1, x_2, x_3) = x_1^4 x_2^3$ and
\[q(x_1, x_2, x_3) = 2x_1^2+3x_2^2+4x_3^2 -1.
\]
Then the following function is harmonic on $\mathbf{R}^3$ and agrees with $p$ on the ellipsoid $\{x \in \mathbf{R}^3: q(x) = 0 \}$:
\[
\begin{split}
x_1^4 &x_2^3 + (2x_1^2 + 3x_2^2 +  4x_3^2 - 1)
   \Bigl( \frac{97950}{20144813} x_2^5- 
     \frac{2524856930}{100139865423} \,x_1^2 x_2\\[8pt]
     &- 
     \frac{3423451}{60434439} \,x_1^4 x_2- 
     \frac{148091}{33379955141}\, x_2^3 - 
     \frac{2306686}{20144813} \,x_1^2 x_2^3- 
     \frac{701980831}{500699327115}\,x_2\\[8pt]
     &+ 
     \frac{32326712}{7703066571}\,x_2 x_3^2 + 
     \frac{3712712}{60434439}\,x_1^2 x_2 x_3^2 + 
     \frac{53836}{20144813} \,x_2^3 x_3^2- 
     \frac{236464}{60434439} \,x_2 x_3^4 \Bigr).\\[8pt]
\end{split}
\]
\end{example}

The function above obviously equals $x_1^4 x_2^3$ on the ellipsoid in question. Thus to verify that the function above is indeed the solution to our Dirichlet problem, it is only necessary to verify that the Laplacian of the function above equals~$0$; \textit{Mathematica} can easily perform this calculation.

Note that, as expected from our discussion in the previous section, every exponent of each $x_j$ in the solution above has the same parity as in the boundary function. Specifically, in the solution above the exponents of $x_1$, $x_3$, and $x_4$ are even and the exponent of $x_2$ is odd, thus following the pattern of the boundary function $x_1^4 x_2^3$.

Our \textit{Mathematica} implementation of the algorithm can handle quadratic surfaces defined with symbols as well as concrete numbers. This capability is illustrated in our second example, which takes place in dimension $4$ with a boundary function of degree~$7$. 

\begin{example} \label{ex2}
Suppose $p(x_1,x_2, x_3, x_4) = x_1^3x_2^2x_3x_4$ and
\[
q(x_1,x_2, x_3, x_4) =
cx_1^2+3x_2^2+4x_3^2+5x_4^2 - 1.
\]
Then the following function is harmonic on $\mathbf{R}^4$ and agrees with $p$ on the ellipsoid $\{x \in \mathbf{R}^4: q(x) = 0 \}$:
\[
\begin{split}
x_1^3&x_2^2x_3x_4 \\[8pt]
  & + \frac{cx_1^2 + 3x_2^2 + 4x_3^2 + 5x_4^2 - 1}
  {788400 + 367920c + 48712c^2 + 2520c^3 + 
          45c^4}
  \Bigl(
  4( 50 + 3c) ( 36 + 5c)   {x_1}x_3^3x_4 \\[8pt]
          &- 
  12( 2190 + 281c + 9c^2 )x_1^3{x_3}{x_4}-\frac{( 82800 + 21868c + 1764c^2 + 
         45c^3)}{3( 10 + c ) } \,{x_1}{x_3}{x_4}  \\[8pt]
  &- ( 50400 + 21118c + 1845c^2 + 
     45c^3 ){x_1}x_2^2{x_3}{x_4}  + 
  5( 46 + 3c ) ( 36 + 5c ) 
   {x_1}{x_3}x_4^3 \Bigr). \\[8pt]
\end{split}
\]
\end{example}

We now return to dimension $3$ but increase the degree of the boundary function to~$10$.

\begin{example} \label{ex3}
The harmonic polynomial that equals $x_1^{10}$ on the ellipsoid
\[
\{x \in \mathbf{R}^3 : 2x_1^2 + 3x_2^2 + 4x_3^2 - 1 = 0\}
\]
has value
\[
\frac{500945213823452554440546462385400584789}
{397263369506735959801289842040922215251461}
\]
at the origin.
\end{example}

Of course, the value at the origin of the solution to this Dirichlet problem is simply the constant term in the polynomial that gives the solution. In the example above, we have given only the constant term because displaying the entire solution would require more than a page.

The example above illustrates the difference between the sphere and ellipsoids. The solution to the Dirichlet problem on the unit sphere in $\mathbf{R}^3$ with boundary function $x_1^{10}$ has value $1/11$ at the origin, and the largest integer appearing in the numerator or denominator of the coefficients of any term of the solution is $46189$. In contrast, if the ellipsoid $\{x \in \mathbf{R}^3 : 2x_1^2 + 3x_2^2 + 4x_3^2 - 1 = 0\}$ replaces the unit sphere, then the value of the solution at the origin is given by the fraction in the example above, and the coefficients of the other terms of the solution have similarly huge numerators and denominators.

We conclude this section by giving some data about the speed of our algorithm. The times given below are CPU times used by the \textit{Mathematica} implementation of our software on a Windows desktop computer with a 1.7~GHz Intel Pentium~4 chip.

Examples \ref{ex1} and \ref{ex2} each took about 0.04 seconds. To produce Example~\ref{ex3}, we found the solution to the Dirichlet problem with boundary function $x_1^{10}$ on the ellipsoid $\{x \in \mathbf{R}^3 : 2x_1^2 + 3x_2^2 + 4x_3^2 - 1 = 0\}$; this took less than 0.2~seconds. On the same ellipsoid, solving the Dirichlet problem with boundary function $x_1^{20}$ takes 55~seconds. Changing the boundary function to $x_1^{25}$ increases the solution time to 10~minutes, 5~seconds. Finally, changing the boundary function to $x_1^{30}$ (still on the same ellipsoid) increases the solution time to 2 hours, 29~minutes, 32~seconds.

The long solution times reported in the paragraph above arise because of the huge overhead associated with exact manipulation of the gigantic numerators and denominators  that appear in the solutions to these high-degree problems. For example, most of the numerators and denominators are larger than $10^{200}$ for the Dirichlet problem above with boundary value $x_1^{20}$, rising to over $10^{450}$ for the Dirichlet problem above with boundary value $x_1^{25}$.

However, our algorithm is fast even for high-degree boundary functions if we use \textit{Mathematica}'s floating point arithmetic instead of exact rational arithmetic. For example, if we ask the \textit{Mathematica} implementation of our algorithm to use floating point arithmetic to solve the Dirichlet problem mentioned above with boundary function $x_1^{30}$, then the computation time is reduced from almost 2.5~hours to find the exact solution to just 8.9~seconds to find a very good decimal approximation of the solution. 

\appendix
\section*{Appendix}
\setcounter{section}{1}
\setcounter{equation}{0}

This purpose of this appendix is to prove the differentiation formula given by Proposition~\ref{theorem:1}, which was used in the derivation of our algorithm in Section~\ref{S:algorithm}.

We begin with following lemma, which gives a formula for differentiating the product of a linear function and a polynomial.

\begin{lemma} \label{L:diffxg}
   Suppose $f$ is a polynomial on $\mathbf{R}^n$ and
    $\alpha$ is a multi-index.  Then
\[
                  D^{\alpha} (g f) = 
                   g D^{\alpha} f+
    \sum _{j=1}^n \alpha_j (D_j g)( D^{\alpha - e_j}f )               
\]
for every $g \in \mathcal{P}_1$.
\end{lemma}

\begin{proof}  Fix $g \in \mathcal{P}_1$. We will prove our desired result by induction on $|\alpha|$. To get started, note that the desired result is obviously true when $|\alpha| = 0$.

Now suppose that $|\alpha| > 0$ and that the desired result holds for all multi-indices of smaller order. Choose $k$ such that $\alpha_k > 0$. Then
\begin{align*}
D^\alpha(gf) &=  D^{\alpha - e_k}D_k(gf) \\
&= D^{\alpha - e_k} \bigl(g D_k f + (D_kg)f\bigr) \\
&= D^{\alpha - e_k} (g D_k f) + (D_kg)(D^{\alpha - e_k} f),
\\
\intertext{where the last equation holds because $D_k g$ is constant. Applying our induction hypothesis to evaluate $D^{\alpha - e_k} (g D_k f)$ now gives}
D^\alpha(gf) &= g D^\alpha f + \sum_{j=1}^n \alpha_j (D_j g)(D^{\alpha - e_j}f),
\end{align*}
completing the proof.
\end{proof}

The next proposition, which gives a formula for differentiating the product of a quadratic polynomial and another polynomial, was used in deriving the system of equations~(\ref{sys}).

\begin{proposition} \label{theorem:1}
Suppose $f$ is a polynomial on $\mathbf{R}^n$ and $\alpha$ is a multi-index. Then
\[
D^\alpha ( q f) =qD^\alpha f + \sum _{j=1}^n \alpha_j (D_j q) (D^{\alpha-e_j} f)
+ \frac{1}{2}\sum _{j=1}^n
\alpha_j(\alpha_j - 1) (D_j^2 q) (D^{\alpha - 2e_j} f)
\]
for every nonhyperbolic quadratic $q$.
\end{proposition}

\begin{proof} Fix a nonhyperbolic quadratic $q$. We will prove our desired result by induction on $|\alpha|$. To get started, note that the desired result is obviously true when $|\alpha| = 0$.

Now suppose that $|\alpha| > 0$ and that the desired result holds for all multi-indices of smaller order. Choose $k$ such that $\alpha_k > 0$. Then
      \begin{align}
           D^{\alpha}(qf)  &= \notag
           D^{\alpha - e_k}D_k(qf) \\
&  =  D^{\alpha - e_k}(qD_k f) +
D^{\alpha - e_k}\bigl((D_kq)f\bigr). \label{app6}
\end{align}
Using Lemma~\ref{L:diffxg} to evaluate the last term (applicable because $D_kq \in \mathcal{P}_1$) gives
\begin{equation} \label{app7}
D^{\alpha - e_k}\bigl((D_k q)f\bigr) =
(D_k q)(D^{\alpha - e_k}f) + (\alpha_k -1)(D_k^2q )(D^{\alpha - 2 e_k}f),
\end{equation}
where we have used the fact that $D_j D_k q = 0$ whenever $j \ne k$.
Using our induction hypothesis to evaluate the first term on the right side of (\ref{app6}) and using (\ref{app7}) to evaluate the second term on the right side of (\ref{app6}) now gives
\[
D^\alpha ( q f) =qD^\alpha f + \sum _{j=1}^n \alpha_j (D_j q) (D^{\alpha-e_j} f)
+ \frac{1}{2}\sum _{j=1}^n
\alpha_j(\alpha_j - 1) (D_j^2 q) (D^{\alpha - 2e_j} f),
\]
completing the proof.
\end{proof}

\end{document}